\documentclass[11pt, preprint]{elsarticle}
\usepackage{amsmath, amscd}
\usepackage{amsfonts}
\usepackage{amssymb}
\usepackage{amsthm}
\usepackage{eufrak}
\usepackage{verbatim}
\usepackage[mathscr]{euscript}

\usepackage{xspace}
\def\ul{\underline}

\def\ds{\displaystyle}
\newcommand{\R}{\ensuremath{\mathbb{R}}\xspace}

\newcommand{\T}{\ensuremath{\mathbb{T}}\xspace}
\newcommand{\Pf}{{\bigskip\noindent\bf Proof.\quad} }

\newcommand{\wtilde}{\widetilde}

\newtheorem{theorem}{Theorem}[section] % numbered like the section
\newtheorem{lemma}[theorem]{Lemma} % numbered like the theorems
\newtheorem{proposition}[theorem]{Proposition}

\theoremstyle{definition} % styled differently... not italicized
\newtheorem{definition}[theorem]{Definition}

\begin{document}
\title{On a 1D transport equation with nonlocal velocity and supercritical dissipation }
\author{Tam Do}
\ead{tamdo@math.wisc.edu}
\address{Dept. of Mathematics, University of Wisconsin, Madison, USA}
\date{\today}
\begin{abstract}
We study a 1D transport equation with nonlocal velocity. First, we prove eventual regularization of the viscous regularization when dissipation is in the supercritical range with non-negative initial data. Next, we will prove global regularity for solutions when dissipation is slightly supercritical. Both results utilize a nonlocal maximum principle.
\end{abstract}

\begin{keyword}
Active scalar \sep Supercritical dissipation \sep Nonlocal velocity \sep Modulus of continuity
\end{keyword}
\maketitle

\section{Introduction}
We will consider solutions to the following initial value problem
\begin{eqnarray}
\label{eq:PDE}
\theta_t = (H\theta)\theta_x-(-\Delta)^\alpha \theta \\
\theta(x,0) = \theta_0(x)  \nonumber
\end{eqnarray}
for $(x,t)\in \R\times [0,\infty)$ and $0\le \alpha <1/2$. Here, $H$ is the Hilbert transform
$$
H\theta(x)=P.V. \frac{1}{\pi}\int \frac{\theta(y)}{x-y}\, dy
$$
The equation has an $L^\infty$ maximum principle which makes the range $0\le \alpha<1/2$ supercritical.
The equation can be thought of as a 1D model for the  2D surface quasi-geostrophic (SQG) equation
\begin{eqnarray*}
\theta_t &=& u\cdot \nabla \theta-(-\Delta)^\alpha \theta \quad\mbox{on}\,\,\, \R^2\times[0,\infty) \\
u &=& \nabla^\perp \Delta^{-1/2} \theta.
\end{eqnarray*}
In addition, $\eqref{eq:PDE}$ also has similarities with the Birkhoff-Rott equations for the evolution of vortex sheets (see \cite{Blow2} for more references). For $\alpha \ge 1/2$, the problem $\eqref{eq:PDE}$ is globally well-posed for initial data in $H^{3/2-2\alpha}$ and locally well-posed for $0\le \alpha <1/2$. \cite{Dong1}. For $0\le \alpha <1/4$, finite time blow up has been shown to be possible \cite{Blow1, Blow2}. 

\medskip
In the range $1/4 \le \alpha <1/2$, it is unknown whether there is finite time blow-up or global regularity. One can think of the term $(-\Delta)^\alpha\theta$ as smoothing while the nonlinear term $(H\theta)\theta_x$ as introducing singular behavior. In the supercritical range, scaling favors the nonlinear term and standard methods for proving regularity are not sufficient. For $1/4\le \alpha < 1/2$, if one were to prove the more interesting result of global regularity, one would need to discover some mechanism of nonlinear depletion present in the equation.

\bigskip
In this paper, we will show that two results that are true for the SQG equation also hold for this 1D model. First, we will show eventual regularization for dissipation in the supercritical regime $0< \alpha <1/2$ with non-negative initial data. Second, we will prove global regularity for the slightly supercritical version of this equation. For the SQG equation, the arguments rely on dissipation in the "perpendicular" direction \cite{EvReg} or incompressibility of the fluid velocity \cite{Dabk, Luis}, which are both absent in our 1D setup. In our results, we will need to carefully use the structure of the nonlinearity as well as the exact formula for the dissipation term. The new ingredients in our proof are Lemma 2.7 and a part of section 3.3.

\section{Eventual Regularization}

In this section, we will closely follow the arguments of Kiselev \cite{EvReg}. We will work with solutions to the  dissipative regularization of $\eqref{eq:PDE}$:
\begin{eqnarray}
\label{eq:reg}
\theta_t = (H\theta)\theta_x-(-\Delta)^{\alpha} \theta+\epsilon \Delta \theta \\
\theta(x,0)=\theta_0(x). \nonumber
\end{eqnarray}
The solutions of $\eqref{eq:reg}$ will be smooth and we will estimate the Holder norms of these solutions uniformly in $\epsilon>0$. The limit obtained by letting $\epsilon\to 0$ will yield a candidate for a weak solution $\theta(x,t)\in C_w([0,\infty);L^2) \cap L^2([0,\infty); H^{1/2})\cap L^\infty([0,\infty); L^2)$. However, this regularity appears to 
be insufficient to conclude that $\theta $ solves $\eqref{eq:PDE}$ in the standard weak sense. The equation is not conservative. On the other hand, the limit will inherit our estimates on the regularization. By having control of high enough Holder norms, the following theorem allows us to conclude smoothness :

\begin{theorem}
Let $\theta$ be a solution of $\eqref{eq:reg}$ with non-negative initial data. Let $\beta>1-2\alpha$ and let $0<t_0<t<\infty$. If $\theta\in L^\infty([t_0,t]; C^\beta(\R))$ then $\theta\in C^\infty((t_0,t]\times \R)$ with bounds independent of $\epsilon$.
\end{theorem}

The proof of Theorem 2.1 is analogous to the proof of Theorem 3.1 in Constantin and Wu \cite{Constantin} where they showed a similar result is true for the SQG equation. Their argument for SQG uses Besov space techniques and does not rely heavily on incompressibility, the key difference between $\eqref{eq:PDE}$ and SQG. Since we have non-negative initial data, solutions of $\eqref{eq:PDE}$ are bounded in $L^2$ without the need of incompressibility \cite{Blow2}, a condition necessary to show the analogous result for SQG. Also, $\eqref{eq:PDE}$ posess the same scaling as SQG and the Hilbert transform is bounded on the Holder and Besov spaces, like the Riesz transforms, which is needed in the proof.

\medskip
The main result of this section is the following theorem.

\begin{theorem}
\label{reg}
 Let $\theta(x,t)$ be the limit obtained in $ \eqref{eq:reg}$ by letting $\epsilon\to 0$ with $\theta_0\in L^\infty\cap H^{3/2-2\alpha}$ and non-negative. Then there exist times $0<T_1(\alpha, \theta_0)\le T_2(\alpha, \theta_0)$ such that $\theta$ is smooth for $0<t<T_1$ and $t>T_2$ ($\theta$ will be a classical 
 solution of $\eqref{eq:PDE}$ for these times). 
 \end{theorem}

\bigskip
\noindent
{\bf Remark.} For $T_1\le t\le T_2$, it is unclear in what sense the $\theta$ above is a solution of $\eqref{eq:PDE}$. The theorem follows from uniform in $\epsilon$ estimates for (2) and such estimates can be regarded as the main result of this section.

\bigskip
To control Holder norms, we will show that a certain family of moduli of continuity is eventually preserved under the evolution. 

\begin{definition} 
A function $\omega(\xi): (0,\infty)\mapsto (0,\infty)$ is a modulus of continuity if $\omega$ is increasing, continuous on $(0,\infty)$, concave, and piecewise $C^2$ with one-sided derivatives defined at every point in $[0,\infty)$. A function $f(x)$ obeys $\omega$ if $|f(x)-f(y)|<\omega(|x-y|)$ for all $x\ne y$.
\end{definition}

\bigskip
To prove that solutions preserve a modulus of continuity we state the following lemma, which describes the scenario in which the modulus is broken.

\medskip
\begin{lemma}
\label{break}
 Let $\theta(x,t)$ be a solution of $\eqref{eq:PDE}$. Suppose that $\omega(\xi,t)$ is continuous on $(0,\infty)\times [0,T]$, piecewise $C^1$ in the time variable (with one-sided derivatives defined at all points) for each fixed $\xi>0$, and that for each fixed $t\ge 0$, $\omega(\xi,t)$ is a modulus of continuity. Assume in addition that for each $t\ge 0$, either $\omega(0+,t)>0$, or $\partial_\xi \omega(0+,t)=\infty$, or $\partial_{\xi\xi}^2\omega(0+,t)=-\infty$, and that $\omega(0+,t)$, $\partial_\xi \omega(0+,t)$ are continuous in $t$ with values in $\R\cup \infty$. Let the initial data $\theta_0(x)$ obey $\omega(\xi,0)$. Suppose that for some $t>0$ the solution $\theta(x,t)$ no longer obeys $\omega(\xi,t)$. Then there exist $t_1>0$ and $x,y\in \R$, $x\ne y$ such that for all $t<t_1$, $\theta(x,t)$ obeys $\omega(\xi,t)$ while
$$
\theta(x,t_1)-\theta(y,t_1)= \omega(|x-y|,t_1).
$$
\end{lemma}

\bigskip
The proof of preceding lemma can be found in \cite{EvReg} for the periodic case. Decay results for solutions from \cite{Dong1} allow the lemma to be extended to the non-periodic setting \cite{Du}. We will use the same moduli of continuity as in \cite{EvReg}:
$$
\omega(\xi,\xi_0) = \left\lbrace \begin{array}{ll} \beta H\delta^{-\beta}\xi_0^{\beta-1}\xi+(1-\beta)H\delta^{-\beta}\xi_0^\beta, & 0<\xi<\xi_0 \\ H(\xi/\delta)^\beta, & \xi_0 \le \xi\le \delta \\
H, & \xi>\delta \end{array}\right.
$$
where $\beta>1-2\alpha$.
Observe that if $2\|\theta_0\|_{L^\infty} \le \omega(0,\delta)=(1-\beta) H$, then $\theta_0$ obeys $\omega(\xi,\delta)$. Thus, for every bounded initial data, we can find a modulus that is obeyed.

\medskip
It is known that for $0<\alpha<1/2$, 
$$
(-\Delta)^\alpha \theta(x) = P.V.\int_{-\infty}^\infty \frac{\theta(x)-\theta(x+y)}{|y|^{1+2\alpha}}\, dy,
$$
see \cite{Cordoba} for a proof.
\noindent
We will need the following estimate of the dissipation terms:
\bigskip
\begin{lemma}
\label{dis}
 (Dissipation Estimate) Let $\xi=|x-y|$. Then
$$
-(-\Delta)^\alpha \theta(x,t) + (-\Delta)^\alpha \theta(y,t) \le D_\alpha(\xi, t)
$$
where
$$
D_\alpha(\xi,t) = c_\alpha \left(\int_0^{\frac{\xi}{2}} \frac{\omega(\xi+2\eta,t)+\omega(\xi-2\eta,t)-2\omega(\xi,t)}{\eta^{1+2\alpha}}\, d\eta + \int_{\frac{\xi}{2}}^\infty \frac{\omega(\xi+2\eta, t)-\omega(2\eta -\xi,t)-2\omega(\xi,t)}{\eta^{1+2\alpha}}\, d\eta\right).
$$
\end{lemma}

\medskip
\noindent
See \cite{EvReg} for the proof. Theorem \ref{reg} is an easy consequence of the following lemma, which we will prove later.

\bigskip
\bigskip

\begin{lemma}
\label{lem}
Assume that $\theta_0(x)$ of $\eqref{eq:PDE}$ obeys $\omega(\xi,\delta)$. Then there exist positive constants $C_{1,2}= C_{1,2}(\alpha,\beta)$ such that if $\xi_0(t)$ is a solution of
\begin{equation}
\label{eq:ODE}
\frac{d\xi_0}{dt}=-C_2 \xi_0^{1-2\alpha}, \quad \xi_0(0)=\delta,
\end{equation}
and $H\le C_1 \delta^{1-2\alpha}$, then the solution $\theta(x,t)$ obeys $\omega(\xi,\xi_0(t))$ for all $t$ such that $\xi_0(t)\ge 0$.
\end{lemma}

\bigskip
\noindent
{\bf Proof of Theorem \ref{reg}} The solution $\xi_0(t)$ of $\eqref{eq:ODE}$ becomes zero and stays zero in finite time. Then eventually, the solution $\theta(x,t)$ obeys $\omega(\xi,0)$ and we can uniformly bound its $C^\beta$ norm, $\beta>1-2\alpha$. $\Box$

\bigskip
To prove lemma \ref{lem}, we will show that the breakthrough scenario described in lemma \ref{break} cannot happen. Suppose there exists $t_1>0$ such that $\theta(x,t)$ obeys $\omega(\xi,t)$ for $t<\xi_0(t_1)$ and $\theta(x,t_1)-\theta(y,t_1)=\omega(\xi,\xi_0(t_1))$ where $\xi=|x-y|$. Then it is not hard to see that $\nabla\theta(x,t_1)=\partial_\xi \omega(\xi,\xi_0(t_1)) \frac{x-y}{\xi}= \nabla\theta(y,t_1)$ and $\Delta \theta(x,t_1)-\Delta\theta(y,t_1) \le 2 \partial_{\xi\xi}^2\omega(\xi,t_1)$ (details are in \cite{EvReg}). Also, by \eqref{eq:reg}, lemma \ref{dis},

\bigskip
\begin{equation}
\label{eq:break}
\partial_t \left. \left[ \frac{\theta(x,t)-\theta(y,t)}{\omega(\xi,\xi_0(t))}\right] \right|_{t=t_1} \le
\frac{\Omega(x,y,t_1)\partial_{\xi}\omega(\xi,\xi_0(t_1)) +d_\alpha(\xi,t_1)+2\epsilon \partial_{\xi\xi}^2 \omega(\xi,t_1) - \partial_t \omega(\xi,\xi_0(t_1))}{\omega(\xi, \xi_0(t_1))}
\end{equation}

\begin{comment}
% old bound
\bigskip
\begin{equation}
\label{eq:break}
\partial_t \left. \left[ \frac{\theta(x,t)-\theta(y,t)}{\omega(\xi,\xi_0(t))}\right] \right|_{t=t_1} \le \frac{\Omega(x,y,t_1)\partial_{\xi}\omega(\xi,\xi_0(t_1)) +D_\alpha(\xi,t_1)+2\epsilon \partial_{\xi\xi}^2 \omega(\xi,t_1) - \partial_t \omega(\xi,\xi_0(t_1))}{\omega(\xi, \xi_0(t_1))}
\end{equation}
\end{comment}

\noindent
where $\Omega(x,y,t_1)=H\theta(x,t_1)-H\theta(y,t_1)$ and
$$
d_\alpha(\xi,t_1)= \frac{1}{2} \left(-(-\Delta)^\alpha \theta(x,t_1)+(-\Delta)^\alpha \theta(y,t_1)\right) + \frac{1}{2} D_\alpha(\xi,t_1).
$$
If we can show that the numerator of the right hand side of $\eqref{eq:break}$ is negative, then the modulus of continuity must have been broken at an earlier time, a contradiction. Because of the concavity of $\omega$, $2\epsilon \partial_{\xi\xi}^2 \omega(\xi,t_1) \le 0$ and since we want our estimates to be independent of $\epsilon$, we will ignore this term.

\begin{lemma}
\label{key}
$$
\Omega(x,y)=H\theta(x)-H\theta(y) \le C\left[ \xi^{2\alpha} \left((-\Delta)^\alpha \theta(x)-(-\Delta)^\alpha \theta(y)\right)+ \xi \int_{\xi/2}^\infty \frac{\omega(r)}{r^2}\, dr\right]
$$
\end{lemma}

For simplicity of expression, we have omitted time in our expressions.

\medskip
\Pf 

The first term on the right side will control the singular behavior of the Hilbert transforms near the kernel singularity and the second term will control the behavior away from the singularity.  Where appropriate, integrals will be understood in the principal value sense. 

\medskip
Let $\wtilde{x}= \frac{x+y}{2}$.  Then
\begin{eqnarray*}
\left|\int_{|x-z|\ge \xi} \frac{\theta(z)}{x-z}\, dz-\int_{|y-z|\ge \xi} \frac{\theta(z)}{y-z}\, dz\right| &=& \left|\int_{|x-z|\ge \xi} \frac{\theta(z)-\theta(\wtilde{x})}{x-z}\, dz-\int_{|y-z|\ge \xi} \frac{\theta(z)-\theta(\wtilde{x})}{y-z}\, dz \right| \\ &\le & \int_{|\wtilde{x}-z|\ge \xi/2} \left|\frac{1}{x-z}-\frac{1}{y-z}\right| |\theta(z)-\theta(\wtilde{x})|\, dz \\ & \le & C\xi \int_{|\wtilde{x}-z|\ge \xi/2} \frac{1}{|\wtilde{x}-z|^2} |\theta(z)-\theta(\wtilde{x})|\, dz \le C\xi \int_{\xi/2}^\infty \frac{\omega(r)}{r^2}\, dr
\end{eqnarray*}

Now, we will estimate the other part. Observe that
\begin{eqnarray*}
H\theta(x) -H\theta(y)= \int_{-\infty}^\infty \frac{\theta(x+z)-\theta(x)}{z}\, dz - \int_{-\infty}^\infty \frac{\theta(y+z)-\theta(y)}{z}\, dz
\end{eqnarray*}

\medskip
Then
\medskip

\begin{eqnarray*}
&\, &\int_{|z|<\xi} \frac{\theta(x+z)-\theta(x)}{z}\, dz - \int_{|z|<\xi} \frac{\theta(y+z)-\theta(y)}{z}\, dz - \xi^{2\alpha}((-\Delta)^\alpha\theta(x)- (-\Delta)^{\alpha} \theta(y)) \\ &=& \int_{|z|<\xi} \frac{\theta(x+z)-\theta(x)}{z}\, dz - \int_{|z|<\xi} \frac{\theta(y+z)-\theta(y)}{z}\, dz - \xi^{2\alpha}\int_{-\infty}^\infty \frac{\theta(x)-\theta(x+z)}{|z|^{1+2\alpha}}\, dz \\ &\, & +\xi^{2\alpha}\int_{-\infty}^\infty \frac{\theta(y)-\theta(y+z)}{|z|^{1+2\alpha}}\, dz \\ &=&
\int_{|z|<\xi} \left(\theta(x+z)-\theta(y+z)+\theta(y)-\theta(x)\right) \left(\frac{1}{z}+\frac{\xi^{2\alpha}}{|z|^{1+2\alpha}}\right) \\
& \, & + \int_{|z|>\xi} \left(\theta(x+z)-\theta(y+z)+\theta(y)-\theta(x)\right) \frac{\xi^{2\alpha}}{|z|^{1+2\alpha}}\, dz \\
&=&\int_{|z|<\xi} \left(\theta(x+z)-\theta(y+z)-\omega(\xi)\right) \left(\frac{1}{z}+\frac{\xi^{2\alpha}}{|z|^{1+2\alpha}}\right)\\ &\,& +\int_{|z|>\xi} \left(\theta(x+z)-\theta(y+z)-\omega(\xi)\right) \frac{\xi^{2\alpha}}{|z|^{1+2\alpha}}\, dz \le 0
\end{eqnarray*}

The last inequality follows from the facts that
$$
\theta(x+z)-\theta(y+z)-\omega(\xi) \le 0
$$
and that in our region of integration
$$
\frac{1}{z}+\frac{\xi^{2\alpha}}{|z|^{1+2\alpha}} \ge \frac{1}{z}+\frac{1}{|z|} \ge 0
$$
Thus, we have control over the Hilbert transforms near the kernel singularity. Combining our estimates, we get the result. $\Box$

\begin{comment}
From this, we will show $\ds \partial_t \left[ \frac{\theta(x,t)-\theta(y,t)}{\omega(\xi,t)}\right] |_{t=t_1}$ is negative ($\xi:=|x-y|$), which is contradiction to the minimality of $t_1$. From the equation,
\begin{equation}
\label{eq:1}
\partial_t(\theta(x,t)-\theta(y,t))|_{t=t_1} = (H\theta(x))\theta_x(x,t_1)-(H(\theta(y))\theta_x(y,t_1)- (-\Delta)^\alpha \theta(x,t_1)+(-\Delta)^\alpha \theta(y,t_1).
\end{equation}
It suffices to show that the above quantity is negative. It is not hard to see that since we are at the time $t_1$ that the modulus of continuity is broken, $\theta_x(x,t_1)=\partial_\xi \omega(\xi,t) = \theta_x(y,t_1)$. Thus, we can rewrite $\eqref{eq:1}$ as
$$
(H\theta(x)-H\theta(y)) \partial_\xi \omega(\xi,t)- (-\Delta)^\alpha \theta(x,t_1)+(-\Delta)^\alpha \theta(y,t_1)
$$
\end{comment}

\bigskip
\bigskip
\noindent
{\bf Proof of Lemma \ref{lem}}

We want to show
\begin{equation}
\label{eq:1}
\partial_t \omega(\xi,t_1) > (H\theta(x,t_1)-H\theta(y,t_1))\partial_{\xi}\omega(\xi,\xi_0(t_1)) +d_\alpha(\xi,t_1)
\end{equation}

\medskip
From Lemma 3.3 of \cite{EvReg}, we can choose the constant $C_2$ in $\eqref{eq:ODE}$ small enough so that we have $\partial_t \omega(\xi,\xi_0(t))> \frac{1}{4}D_\alpha (\xi,t)$ at $t=t_1$ ($\xi_0'(t)$ is small). 
By Lemma 5.3 of \cite{EvReg}, we can replace $\ds \xi \int_{\xi/2}^\infty \frac{\omega(r)}{r^2}\, dr$ in Lemma 2.7 by $\omega(\xi,\xi_0)$. Using an argument very similar to Lemma 3.3 of \cite{EvReg}, it can be shown that for all $0<\xi<\delta$,
$$
C\omega(\xi,\xi_0(t_1))\partial_\xi \omega(\xi,\xi_0(t_1)) \le - \frac{1}{4} D_\alpha (\xi,t_1).
$$
where $C$ is the constant from Lemma 2.7.
Now, for $0<\xi<\delta$, we have
\begin{eqnarray*}
C\xi^{2\alpha}  \partial_{\xi}\omega(\xi,\xi_0(t_1)) &\le & C\beta H\xi^{2\alpha} \delta^{-\beta}\xi^{\beta-1} = C\beta \frac{H}{\delta^{1-2\alpha}} \left(\frac{\xi}{\delta}\right)^{2\alpha+\beta-1} 
\end{eqnarray*}

% \left((-\Delta)^\alpha \theta(x,t_1)-(-\Delta)^\alpha \theta(y,t_1)\right)

\begin{comment}

%bad estimate

\begin{eqnarray*}
C\xi^{2\alpha} \left((-\Delta)^\alpha \theta(x,t_1)-(-\Delta)^\alpha \theta(y,t_1)\right)\partial_\xi \omega(\xi,\xi_0(t_1)) & \le & -C\xi^{2\alpha} D_\alpha(\xi,t_1)\partial_\xi \omega(\xi,\xi_0(t_1))  \\
& \le & -C \beta \frac{H}{\delta^\beta} \xi^{\beta+2\alpha-1} D_\alpha(\xi,t_1) \\
&=& -C\beta \frac{H}{\delta^{1-2\alpha}} \left( \frac{\xi}{\delta}\right)^{\beta+2\alpha-1} D_\alpha(\xi,t_1)
\end{eqnarray*}

\end{comment}

By choosing $C_1$ in $H\le C_1 \delta^{1-2\alpha}$ small enough we can bound the expression above by $\ds \frac{1}{2}$. Then
$$
C\xi^{2\alpha} \partial_{\xi}\omega(\xi,\xi_0(t_1)) \left((-\Delta)^\alpha \theta(x,t_1)-(-\Delta)^\alpha \theta(y,t_1)\right) \le \frac{1}{2} \left((-\Delta)^\alpha \theta(x,t_1)-(-\Delta)^\alpha \theta(y,t_1)\right)
$$
Combining these estimates with Lemma \ref{key}, we have $\eqref{eq:1}$.  $\Box $

\bigskip
\bigskip

\section{Well-posedness for Slightly Supercritical Hilbert Model}

\bigskip 
In this section, we prove global regularity for our model
for which the dissipation can be supercritical by a logarithm. Specifically, we will look at solutions of the
following equation
\begin{equation}
\label{eq:super}
\theta_t= (H\theta)\theta_x-\mathscr{L}\theta , \quad \theta(x,0)=\theta_0(x)
\end{equation}
for $\theta_0\in H^{3/2}(\R)$, where $\mathscr{L} \theta = \frac{(-\Delta)^{1/2}}{\log(1-\Delta)}\theta$ is a Fourier multiplier operator with multiplier 
$$P(\xi)= \frac{|\xi|}{\log(1+|\xi|^2)}.
$$ 
For simplicity, we will only concern ourselves with a dissipative operator of this form. The results of this section can easily be generalized to other similar dissipative operators.
The main result of this section is the following

{\theorem Assume that $\theta_0\in H^{3/2}(\R)$. Then there exists a unique smooth solution $\theta$ of $\eqref{eq:super}$.}

\medskip
First, we have local existence of smooth solutions that we will eventually show can be extended.
\bigskip
{\proposition Let $0<\alpha< 1/2$ and $\theta_0\in H^{3/2}$. Then there exists $T>0$ such that $\eqref{eq:super}$ has a unique solution $\theta$ up to time $T$ that satisfies
$$
\sup_{0<t<T} t^{\beta/(2\alpha)} \|\theta(t,\cdot)\|_{\dot{H}^{3/2-2\alpha+\beta}} <\infty
$$
for any $\beta\ge 0$ and
$$
\lim_{t\to 0} t^{\beta/(2\alpha)} \|\theta(t,\cdot)\|_{\dot{H}^{3/2-2\alpha+\beta}}=0
$$
for any $\beta>0$. Furthermore, we can extend the solution beyond $T$ if $\|\nabla \theta\|_{L^1(0,T;L^\infty)}<\infty$.}

\Pf This result is analogous to Theorem 4.1 and Proposition 6.2 of Dong \cite{Dong1} where it is done for the usual fractional laplacian dissipation. The argument for the dissipation we are using is very similar. We will present the modification necessary to make their proof work. The general idea is that $\mathscr{L}$ is more dissipative then $(-\Delta)^\alpha$ for $0<\alpha<1/2$. Let $\theta$ be a solution of $\eqref{eq:super}$ and let $\theta_j=\Delta_j \theta$ be the $j$th Littlewood-Paley projection. Applying $\Delta_j$ to both sides of $\eqref{eq:super}$ we get
\begin{equation}
\label{paley}
\partial_t \theta_j+(H\theta) (\theta_j)_x+ \mathscr{L}\theta_j = [H\theta,\Delta_j]\theta_x
\end{equation}
where $[H,\Delta_j]$ is a commutator with $ [H\theta,\Delta_j]\theta_x= (H\theta)(\theta_j)_x - \Delta_j((H\theta)(\theta_x))$. By applying Plancherel and using that $\Delta_j$ localizes $\theta$ in the frequency space,
$$
\int_\R \theta_j \mathscr{L}\theta_j\, dx \ge 2^{2\alpha j} C \|\theta_j\|_{L^2}^2
$$
for some constant $C$. Then by multiplying both sides of $\eqref{paley}$ by $\theta_j$ and integrating we get
$$
\frac{1}{2} \frac{d}{dt} \|\theta_j\|_{L^2}^2 + 2^{2\alpha j} C \|\theta_j\|_{L^2}^2 \le \int_\R \left([H\theta, \Delta_j]\theta_x+ H\theta_x \theta_j/2\right)\theta_j\, dx.
$$
This is the same type of inequality used in Dong \cite{Dong1} and one can apply the methods there to arrive at the a priori bounds needed to conclude  local existence as well as higher regularity despite the absence of a divergence free property

\medskip
Dong also proves a Beale-Kato-Majda type blow up criterion for $\eqref{eq:PDE}$ and the result still holds true for $\eqref{eq:super}$ with our logarithmic dissipation. The contribution from the dissipation term is still non-negative, which is the only fact used about dissipation in his proof. Specifically, by Plancherel, for any regular enough function $f$,
$$
\int_\R f\mathscr{L}f\, dx \ge 0
$$
$\Box$

\bigskip
Thus, if we can show $\|\nabla \theta\|_\infty$ is bounded uniformly in time, then Theorem 3.1 is proved. Having such a bound will allow us to extend local solutions indefinitely. To have a bound on $\|\nabla \theta\|_\infty$, we will show the evolution preserves a family of moduli of continuity. If a function $f\in C^2 (\R)$ obeys a modulus $\omega$ satisfying $\omega'(0)<\infty$ and $\omega''(0)=-\infty$, then $\|\nabla f(x)\|_\infty< \omega'(0)$(see \cite{SQG}).  Therefore, if $\theta$ preserves a modulus of continuity, $\|\nabla \theta\|_\infty< \omega'(0)$.

\subsection{Writing $\mathscr{L}$ as dissipative nonlocal operator}

\bigskip
In the proofs, it will be easier to write $\mathscr{L}$ as a nonlocal dissipative nonlocal operator, which the following version of lemmas 5.1 and 5.2 from \cite{Super} allows us to do.

\bigskip
{\lemma The operator $\mathscr{L}$ can be written as
$$
\mathscr{L}\theta(x)=\int_\R (\theta(x)-\theta(x+y))K(y)\, dy.
$$
Also, there exists a positive constant $C$ such that 
$$
\frac{1}{C} \frac{1}{|y|} P(|y|^{-1}) \le K(y) \le C\frac{1}{|y|} P(|y|^{-1})
$$
where the lower bound holds for $|y|<2\sigma$ for some small constant $\sigma$.
}

\medskip
Since are not assured positivity of the kernel $K$, by the previous lemma, we will not have the $L^\infty$ maximum principle. The following result (Lemma 5.4 from \cite{Super}) allows us to circumvent this.

\bigskip
{\lemma Let $\theta$ solve $\eqref{eq:super}$. Then there exists a constant $M=M(P,\theta_0)$ such that $\|\theta(\cdot,t)\|_{L^\infty} \le M$ for all $t\ge 0$.}

\bigskip
Using the notation from Lemma 3.3, let $\varphi$ be a smooth radially decreasing function that is identically $1$ on $|y|\le \sigma$ and vanishes identically on $|y|\ge 2\sigma$.  Let 
\begin{eqnarray*}
K_1(y) &=& K(y)\varphi(y) \\
K_2(y) &=& K(y)(1-\varphi(y))
\end{eqnarray*}
Now, we decompose the dissipation term $\mathscr{L}$:
\begin{equation}
\label{eq:decomp}
\mathscr{L}\theta(x)= \mathscr{L}_1\theta(x)+\mathscr{L}_2\theta(x) := \int_\R (\theta(x)-\theta(x+y))K_1(y)\, dy+ \int_\R (\theta(x)-\theta(x+y))K_2(y)\, dy.
\end{equation}
Let 
$$
m(r)=\frac{1}{C} P(r^{-1}) \varphi(r)
$$
where $C$ is the constant from Lemma 3.3. Then we have the following lower bound on $\mathscr{L}_1$ that we will use extensively:
$$
\mathscr{L}_1\theta(x)\ge  \int_\R (\theta(x)-\theta(x+y))\frac{m(|y|)}{|y|}\, dy
$$
The operator $\mathscr{L}_1$ satisfies the following conditions satisfied by more general nonlocal dissipative operators

\bigskip
\begin{enumerate}
\item there exists a positive constant $C_0>0$ such that 
$$
rm(r)\le C_0\,\,\mbox{for all}\,\, r\in (0,2\sigma)
$$

for some $r_0>0$.

\item there exists some $a>0$ such that $r^{a} m(r)$ is non-increasing.
\end{enumerate}

\bigskip
We also have the following dissipation estimate whose proof is analogous to Lemma \ref{dis}.
{\lemma Suppose $\theta$ obeys a modulus of continuity $\omega$. Suppose there exists $x,y$ with $|x-y|=\xi>0$ such that $\theta(x)-\theta(y)=\omega(\xi)$. Then
$$
\mathscr{L}_1\theta(x)-\mathscr{L}_1\theta(y) \ge \mathscr{D}(\xi)
$$
where
\begin{eqnarray*}
\mathscr{D}(\xi) = &A& \int_0^{\xi/2} \left(2\omega(\xi)-\omega(\xi+2\eta)-\omega(\xi-2\eta)\right) \frac{m(2\eta)}{\eta}\, d\eta \\ &+& A\int_{\xi/2}^\infty \left(2\omega(\xi)-\omega(\xi+2\eta)+\omega(2\eta-\xi)\right) \frac{m(2\eta)}{\eta}\, d\eta
\end{eqnarray*}
and $A$ is a constant.}

\bigskip
\subsection{The Moduli of Continuity}

\bigskip
The modulus from \cite{Super} will work here. Fix a small constant $\kappa >0$. For any $B\ge 1$, define $\delta(B)$ to be the solution of 
$$
m(\delta(B)) = \frac{B}{\kappa}.
$$
We can also assume that $\delta(B)\le \sigma/2$ by choosing $\kappa$ small enough. Let $\omega_B(\xi)$ be a continuous function with $\omega_B(0)=0$ and
\begin{eqnarray}
\omega_B'(\xi) &=& B - \frac{B^2}{2C_a \kappa} \int_0^\xi \frac{3+ \log(\delta(B)/\eta)}{\eta m(\eta)}\, d\eta, \quad \mbox{for}\,\, 0<\xi <\delta(B), \\ 
\omega_B'(\xi) &=& \gamma m(2\xi), \,\,\,\qquad\qquad\qquad\qquad\qquad\qquad \mbox{for}\,\, \xi>\delta(B),
\end{eqnarray}

\medskip
\noindent
where $C_a = (1+3a)/a^2$ and $\gamma>0$ is a constant dependent on $\kappa, A,$ and $m$. It is shown in \cite{Super} that $\omega_B$ is a indeed a modulus of continuity.

\bigskip
Now, we will show that solutions will initially obey some $\omega_B(\xi)$ for some $B$ large enough. Since evolution immediately smooths out the initial data, we can assume $\theta_0$ is a smooth as needed. By Lemma 3.4, it suffices to find $B$ such that $\omega_B(\xi)\ge \min\{\xi \|\nabla \theta_0\|_{L^\infty}, 2M\}$
for all $\xi>0$ where $M$ is from Lemma 3.4. By concavity of $\omega$, we are left to show that $\omega_B(b)\ge 2M$ where $b = 2M/\|\nabla \theta_0\|_{L^\infty}$. Choose $B$ large so $b>\delta(B)$, so
$$
\omega_B(b)= \omega_B(\delta(B))+\int_{\delta(B)}^b \omega_B'(\eta)\, d\eta \ge \gamma \int_{\delta(B)}^b m(2\eta)\, d\eta \to \infty
$$
as $\delta(B)\to 0$. By choosing $B$ possibly even larger we can have $\omega_B(\sigma)\ge 2M\ge 2\|\theta(\cdot,t)\|_{L^\infty}$ where $\sigma$ is from our decomposition of $\mathscr{L}$ earlier. Therefore, the modulus can only be broken for $0<\xi<\sigma$ and solutions will initially obey a modulus from the family $\{\omega_B\}_{B\ge 1}$.

\bigskip
\subsection{The moduli are preserved} 
\medskip
To prove a modulus of continuity is preserved, we will rule out the breakthrough scenario described in Lemma \ref{break}. Let $t_1$ be the time of breakthrough. By using $\eqref{eq:super}$ and Lemma 3.5,
$$
\partial_t\left.(\theta(x,t)-\theta(y,t))\right|_{t=t_1} \le  \left( H\theta(x,t_1)-H\theta(y,t_1)\right) \omega'(\xi) -\mathscr{D}_B(\xi) + \mathscr{L}_2\theta(y,t_1)-\mathscr{L}_2\theta(x,t_1)
$$
where $\xi = |x-y|$.
If the right side of the equation above is negative then the modulus was broken at an earlier time, a contradiction. In \cite{Super}, they show 
$$
|\mathscr{L}_2\theta(x,t)-\mathscr{L}_2\theta(y,t)|\le \frac{1}{2} \mathscr{D}_B(\xi)
$$
for $0<\xi<\sigma$ so to prove Theorem 3.1, it suffices to show
\begin{equation}
\label{eq:show}
\left( H\theta(x,t_1)-H\theta(y,t_1)\right) \omega_B'(\xi) -\frac{1}{2}\mathscr{D}_B(\xi) < 0
\end{equation}
for $0<\xi<\sigma$ where $\mathscr{D}_B$ is the expression from Lemma 3.5 with $\omega_B$ being the modulus. For simplicity, we will now omit $t_1$ from our expressions involving $\theta$.

\bigskip
\noindent
{\bf \ul{Case: $\xi \ge \delta(B)$}} 

\medskip
By a similar argument to the proof of Lemma \ref{key}, 
$$
\left|\int_{|x-z|\ge 2\xi} \frac{\theta(z)}{x-z}\, dz-\int_{|y-z|\ge 2\xi} \frac{\theta(z)}{y-z}\, dz\right|  \le C\xi \int_\xi^\infty \frac{\omega_B(\eta)}{\eta^2}\, d\eta.
$$
Integrating by parts
$$
\xi \int_\xi^\infty \frac{\omega_B(\eta)}{\eta^2}\, d\eta = \omega_B(\xi)+\gamma \xi \int_\xi^\infty \frac{m(2\eta)}{\eta}\, d\eta.
$$
By property (2) of $m$,
$$
\int_\xi^\infty \frac{m(2\eta)}{\eta}\, d\eta \le \xi^a m(2\xi)\int_\xi^\infty \frac{1}{\eta^{1+a}}\, d\eta \le \frac{m(2\xi)}{a}
$$
Now, we have
$$
\xi \int_\xi^\infty \frac{\omega_B(\eta)}{\eta^2}\, d\eta \le \omega_B(\xi) + \frac{\gamma\xi m(2\xi)}{a}
$$
For $\delta(B) \le \xi \le 2\delta(B)$, it is not hard to see that 
$$
\frac{\gamma\xi m(2\xi)}{a} \le \omega_B(\xi),
$$
the details are in \cite{Super}. For $\xi > 2\delta(B)$, we have $\xi -\delta(B) \ge \xi/2$ so
$$
\omega_B(\xi) \ge \gamma \int_{\delta(B)}^\xi m(2\eta)\, d\eta \ge \gamma m(2\xi)(\xi-\delta(B)) \ge \frac{\gamma \xi m(2\xi)}{2}.
$$
Thus, we have,
$$
C\xi \int_\xi^\infty \frac{\omega_B(\eta)}{\eta^2}\, d\eta \le C\left(1+ \frac{2}{a}\right) \omega_B(\xi).
$$
In \cite{Super}, they prove the following estimate on the dissipation term
$$
-\mathscr{D}_B (\xi) \le -\frac{2-c_a}{C} \omega_B(\xi) m(2\xi)
$$
where $c_a=1+(3/2)^{-a}$.
Then we obtain
$$
\left(\int_{|x-z|\ge 2\xi} \frac{\theta(z)}{x-z}\, dz-\int_{|y-z|\ge 2\xi} \frac{\theta(z)}{y-z}\, dz\right) \omega_B'(\xi) - \frac{1}{4} \mathscr{D}_B(\xi) \le \left(C\gamma \frac{a+2}{a} - \frac{2-c_a}{4C}\right) \omega_B(\xi) m(2\xi) < 0
$$
if we set $\gamma$ small enough. In other words, we have used some of the dissipation to control the modulus of the Hilbert transform away from the kernel singularity. Now, we will concern ourselves with the other part of the Hilbert transform. A novel step is that instead of using $\mathscr{D}_B$ we will use the expression for $\mathscr{L}_1\theta$ directly. We want to show
\begin{equation}
\label{eq:show1}
\left(\int_{|x-z|\le 2\xi} \frac{\theta(z)}{x-z}\, dz-\int_{|y-z|\le 2\xi} \frac{\theta(z)}{y-z}\, dz\right) \omega_B'(\xi) - \frac{1}{4}( \mathscr{L}_1\theta(x)-\mathscr{L}_1\theta(y)) < 0.
\end{equation}
After a similar manipulation as in the proof of Lemma 2.6, the left side of $\eqref{eq:show1}$ is precisely
\begin{eqnarray*}
&\, & \int_{|z| <2\xi} \left(\theta(y)-\theta(z+y)-\theta(x)+\theta(x+z)\right) \left[ \frac{\omega_B'(\xi)}{z} +\frac{1}{4} \frac{m(z)}{|z|}\right]\, dz \\
&\, & + \int_{|z|\ge 2\xi}  \left(\theta(y)-\theta(z+y)-\theta(x)+\theta(x+z)\right) \frac{1}{4} \frac{m(z)}{|z|}\, dz
\end{eqnarray*}
By hypothesis, we have $\theta(y)-\theta(x)= \omega_B(\xi)$, so 
$$
\theta(y)-\theta(z+y)-\theta(x)+\theta(x+z)= \omega_B(\xi)- \theta(z+y)+\theta(x+z)< 0.
$$
If we can show 
$$
\frac{\omega_B'(\xi)- \frac{1}{4}m(z)}{z} > 0
$$
for $-2\xi<z<0$ then we are done. However,
$$
\omega_B'(\xi)- \frac{1}{4}m(z) = \gamma m(2\xi)-\frac{1}{4}m(z) < 0
$$
from the fact that $m$ is non-increasing and choosing $\gamma$ small enough. Therefore, we have $\eqref{eq:show}$ and the case when $\xi \ge \delta(B)$ is complete.

\bigskip
\noindent
{\bf \ul{Case: $0< \xi \le \delta(B)$}} The argument for this case is exactly the same as in \cite{Super} with no modifications. Therefore, the proof of theorem 3.1 is complete. $\Box$

\bigskip
\noindent
\ul{{\bf Acknowledgements}} The author wishes to thank Prof. Alex Kiselev for the introduction into the subject and Prof. Hongjie Dong for helpful comments. The author also acknowledges the support of the NSF grant DMS 1147523 and DMS 1159133 at UW-Madison. 

\begin{comment}
\medskip
Since $m$ is singular at the $0$, by choosing $A$ small enough
$$
\omega_B'(\xi)-\frac{1}{2}m(z) \le B-\frac{1}{2}m(z) < 0
$$
for $|z| \le A\xi$. Then 
\begin{equation*}
\left(\int_{|x-z|\le A\xi} \frac{\theta(z)}{x-z}\, dz-\int_{|y-z|\le A\xi} \frac{\theta(z)}{y-z}\, dz\right) \omega_B'(\xi) - \frac{1}{2}( \mathscr{L}\theta(x)-\mathscr{L}\theta(y)) < 0.
\end{equation*}
Also, by an argument similar to the previous case,
$$
\left(\int_{|x-z|\ge A\xi} \frac{\theta(z)}{x-z}\, dz-\int_{|y-z|\ge A\xi} \frac{\theta(z)}{y-z}\, dz\right) \omega_B'(\xi) - \frac{1}{2}( \mathscr{L}\theta(x)-\mathscr{L}\theta(y)) < 0.
$$
Thus, the argument for $0<\xi\le \delta(B)$ is complete.
\end{comment}

\bibliographystyle{plain}
\bibliography{refs}

\end{document}